  \documentclass{ajour}
\usepackage{amsmath}
\usepackage{graphicx}

\def\text{\hbox} \def\newpage{\vfill\break}
\def\z{\psi}
\def\i{\iota}
\def\a{a}
\def\b{b}
\def\d{\delta}

\def\s{\sigma}
\def\S{\Sigma}

\def\te{\theta}
\def\TE{\Theta}
\def\o{\omega}
\def\O{\Omega}
\def\i{\iota}
\def\t{\rho}
\def\v{\upsilon}
\def\g{\gamma}
\def\G{\Gamma}
\def\d{\delta}
\def\D{\Delta}
\def\E{E(\Gamma)}
\def\V{V(\Gamma)}
\def\FF{{\mathcal F}}
\def\GG{{\mathcal G}}
\def\B{{\mathcal B}}
\def\Z{{\mathbf Z}}
\def\xu{\{0,1\}}
\def\xd{(\{0,1\},2)}
\def\xt{(\{0,1\},3)}
\def\pqr{(\{p,q\},r)}
\def\hq{(h_0,h_1,h_2;q_0,q_1,q_2)}
\def\hqp{(h_0',h_1',h_2';q_0',q_1',q_2')}
\def\h{\widehat}
\def\wt{\widetilde}

\begin{document}

%%  Authors, start here ==>>

%\draft % Optional, will cause a line at the bottom of each page
%% with the words `Draft' and the time and date that the article
%% was LaTeXed. Will also double space text.

\title{2-SYMMETRIC TRANSFORMATIONS FOR 3-MANIFOLDS OF GENUS TWO
\thanks{Work performed under the auspices of G.N.S.A.G.A. of
C.N.R. of Italy; supported by the M.U.R.S.T. of Italy, within the
project {\it Topologia e Geometria} and by the University of
Bologna, funds for selected research topics.\\ \indent {\it 1991 Mathematics Subject
Classification:} Primary 57M12, 57M15; Secondary 57M25, 05C10.\\
\indent {\it Keywords:} Manifolds, genus, crystallizations.}
}

\author{Luigi Grasselli}

\affil{Department of Sciences and Methods for Engineering,
University of Modena and Reggio Emilia, 42100 Reggio Emilia,
ITALY}

\email{grasselli.luigi@unimo.it}

\author{Michele Mulazzani\thanks{Also, C.I.R.A.M., Bologna, ITALY.}}

\affil{Department of Mathematics, University of Bologna, 40127 Bologna, ITALY}
\email{mulazza@dm.unibo.it}

\and

\author{Roman Nedela\thanks{Supported in part by the
Ministry of Education of Slovakia,
grant no. 1/3213/96.}}

\affil{Department of Mathematics, M. Bel University, 97549
Bansk\'a Bystrica, SLOVAKIA}

\email{nedela@bb.sanet.sk}

\authorrunninghead{Grasselli, Mulazzani, and Nedela}
\titlerunninghead{2-symmetric transformations}

\abstract{
As previously known, all $3$-manifolds of genus two can be
represented by edge-coloured graphs uniquely defined by $6$-tuples
of integers satisfying simple conditions. The present paper
describes an ``elementary transformation'' on these $6$-tuples
which changes the associated graph but does not change the
represented manifold. This operation is a useful tool in the
classification problem for $3$-manifolds of genus two; in fact, it
allows to define an equivalence relation on ``admissible''
$6$-tuples so that equivalent $6$-tuples represent the same
manifold. Different equivalence classes can represent the same
manifold; however, equivalence classes ``almost always'' contain
infinitely many $6$-tuples. Finally, minimal representatives of
the equivalence classes are described.
}

\begin{article}

%%%%%%%%%%%%%%%%%%%%%%%%%%%%%%%%%%%%%%%%%%%%%%%%%%%%%%%%%

\section{Introduction}

A classical invariant for closed 3-manifolds is the Heegaard genus
\cite{He}. It is well known that the only 3-manifold of genus 0 is the
3-sphere and the 3-manifolds of genus 1 are $\mathbf S^1\times\mathbf
S^2$, $\mathbf S^1\widetilde{\times} \mathbf S^2$ and the lens spaces.
Thus, all 3-manifolds of genus $<2$ are completely classified. On
the contrary, the classification problem for the 3-manifolds of
genus $g$ is still unsolved for $g\ge 2$. The present paper deals
with the classification problem for the class $\mathcal M_2$ of
all orientable 3-manifolds of genus two.

PL-manifolds can be represented by edge-coloured
graphs \cite{FGG} and within this theory the homeomorphism problem between manifolds can be
translated into an
equivalence criterion for edge-colored graphs by means of the so-called
``dipole moves'' \cite{FG}; namely, two manifolds are homeomorphic if and only if
each pair of coloured graphs representing them can be joined by a finite
sequence of dipole moves.

The Heegaard genus of a $3$-manifold can also be
defined in terms of coloured graphs; in fact, the Heegaard genus of a
3-manifold $M$ is the non-negative integer
$g(M) = \text{min}\,\{ g(G) \mid G \text{ represents } M\}$ ,
where $g(G)$ is the minimal genus of a surface into which the coloured graph $G$
``regularly'' embeds (\cite{G1}, \cite{G2}).

In particular, each manifold of $\mathcal M_2$ can be represented
by ``highly symmetric'' graphs, which are uniquely defined by
6-tuples of integers. The classification problem in $\mathcal M_2$
then translates into determining when two $6$-tuples represent the
same manifold. Unfortunately, the dipole moves generally modify
the genus of a coloured graph; hence, single dipole moves cannot
be used for defining an equivalence criterion on 6-tuples, which
translates the homeomorphism of the represented manifolds of
$\mathcal M_2$.

We point out that, up to now, the problem of finding a complete
set of moves translating the homeomorphism between manifolds in
$\mathcal M_2$ is still open in all known representation theories
for $\mathcal M_2$.

Our paper describes an ``elementary transformation'' on 6-tuples
representing the manifolds of $\mathcal M_2$ which changes the
associated graph but does not change the represented manifold;
this is performed by standard sequences of dipole moves which do
not change both the genus and the symmetry of the coloured graph.

This elementary transformation allows us to define an equivalence
relation on $6$-tuples so that equivalent $6$-tuples represent the
same manifold. Different equivalence classes can represent the
same manifold; however, the transformation seems to be a useful
tool for computer generating of reduced catalogues of $\mathcal
M_2$. In fact, we show that ``almost every manifold'' in $\mathcal
M_2$ can be represented by infinitely many equivalent $6$-tuples;
moreover, we describe the minimal representatives of the
equivalence classes.

\section {Preliminaries}

Throughout this paper, all spaces and maps are piecewise-linear (PL)
in the sense of \cite{RS}. Manifolds are always assumed to be closed, connected
and orientable.
For basic graph theory, we refer to \cite{Ha}.
We shall use the term graph instead of multigraph: hence,
loops are forbidden but multiple edges are allowed.

An {\it edge-coloring\/} on a graph $\G=(\V,\E)$ is a map $\g
:\E\to\D_n=\{0,1,\ldots ,n\}$ such that $\g(e)\not=\g(f)$, for
each pair of adjacent edges $e,f$. If $v,w$ are the vertices of an
edge $e\in E(\G)$ such that $\g(e)=c$, we say that $e$ is a {\it
$c$-edge\/} and that $v,w$ are {\it $c$-adjacent\/}. The pair
$(\G,\g)$, $\G$ being a graph and $\g :\E\to\D_n$ being an
edge-coloring, is said to be an $(n+1)$-{\it coloured graph with
boundary\/}. A boundary-vertex is simply a vertex $v$ of degree
less than $n+1$; if there are no $c$-edges incident with $v$, we
say that $v$ is a boundary vertex with respect to colour $c$. If
$\G$ is regular of degree $n+1$ (i.e., if $\G$ has no boundary
vertices), then $(\G,\g)$ is simply called an $(n+1)$-{\it
coloured graph\/}. The notion of {\it colour preserving
isomorphism\/} ({\it c.p.-isomorphism\/}) between $(n+1)$-coloured
graphs is straightforward.

For each $\B\subseteq\D_n$, we set $\G_{\B}=(\V,\g^{-1}(\B))$;
moreover, each connected component of $\G_{\B}$ will be called a
{\it $\B$-residue}. For each colour $c\in\D_n$, we set $\h
c=\D_n-\{c\}$. For sake of conciseness, we shall often denote
$(\G,\g)$ simply by the symbol $\G$ of its underlying graph.

As shown in \cite{FGG}, every $(n+1)$-coloured graph (with boundary) $\G$
represents an $n$-dimensional pseudocomplex $K(\G)$ \cite{HW}, which is
a pseudo-manifold (with boundary) \cite{ST}; moreover, $K(\G)$ is
orientable if and only if $\G$ is bipartite. An {\it $n$-gem\/} is
an $(n+1)$-coloured graph representing an $n$-manifold.

An $(n+1)$-coloured graph $\G$ is said to be {\it contracted\/} if
$\G_{\h c}$ is connected, for each $c\in\D_n$. A {\it
crystallization\/} is a contracted gem. Every $n$-manifold admits
a crystallization \cite{Pe}.

Let $\G$ be an $(n+1)$-coloured graph and let $\te$ be the
subgraph composed by two vertices $X,Y$ joined by $h$ edges ($1\le
h\le n$) with colours $c_1,\ldots,c_h$. If $X$ and $Y$ belong to
distinct components of $\G_{\D_n-\{c_1,\ldots,c_h\}}$, then $\te$
is called a {\it dipole of type $h$\/}.

{\it Cancelling\/} $\te$ means: \begin{itemize}
\item deleting the vertices and the edges of $\te$;
\item welding the ``hanging'' edges of the same colour.\end{itemize}

\noindent {\it Adding\/} $\te$ means the inverse process.
If $\G$ and $\G'$ are $n$-gems of the $n$-manifolds $M$ and $M'$ respectively, then $M'$ is
homeomorphic to $M$ if and only if $\G'$ is obtained from
$\G$ by cancelling and/or adding a finite number of dipoles \cite{FG}.

For a general survey on manifold representation theory by means of
coloured graphs, see \cite{BM}, \cite{FGG}, \cite{Li} and \cite{Vi}.

\section {Blocks and gluing subgraphs}

Let $\G$ be a $4$-coloured graph and let $p,q,r$ be distinct
colours of $\D_3$. Suppose that $C'$ and $C''$ are distinct
$\{p,q\}$-residues of $\G$ and that $v_1',\ldots,v_h'$ (resp.
$v_1'',\ldots,v_h''$) are distinct consecutive vertices of a
$\{p,q\}$-residue $C'$ (resp. $C''$); this means that, for each
$i=1,\ldots,h-1$, $v_i'$ (resp. $v_i''$) is joined with $v_{i+1}'$
(resp. $v_{i+1}''$) by an edge  $\a_i'$ of $C'$ (resp. $\a_i''$ of
$C''$). Moreover, suppose that $\g(\a_i')=\g(\a_i'')$, for each
$i=1,\ldots,h-1$, and $v_j'$ is joined with $v_j''$ by an
$r$-coloured edge $\b_j$, for each $j=1,\ldots,h$. Then, the
subgraph $\O$ of $\G$ defined by:
$$V(\O)=\{v_1',\ldots,v_h',v_1'',\ldots,v_h''\},$$
$$E(\O)=\{\a_1',\ldots,\a_{h-1}',\a_1'',\ldots,\a_{h-1}'',\b_1,\ldots,\b_h\},$$
is called a {\it $\pqr$-block\/} of length $h$, connecting $C'$
with $C''$ (see Figure 1).

\bigskip

\begin{figure}[bht]
 \begin{center}
 \includegraphics*[totalheight=3.5cm]{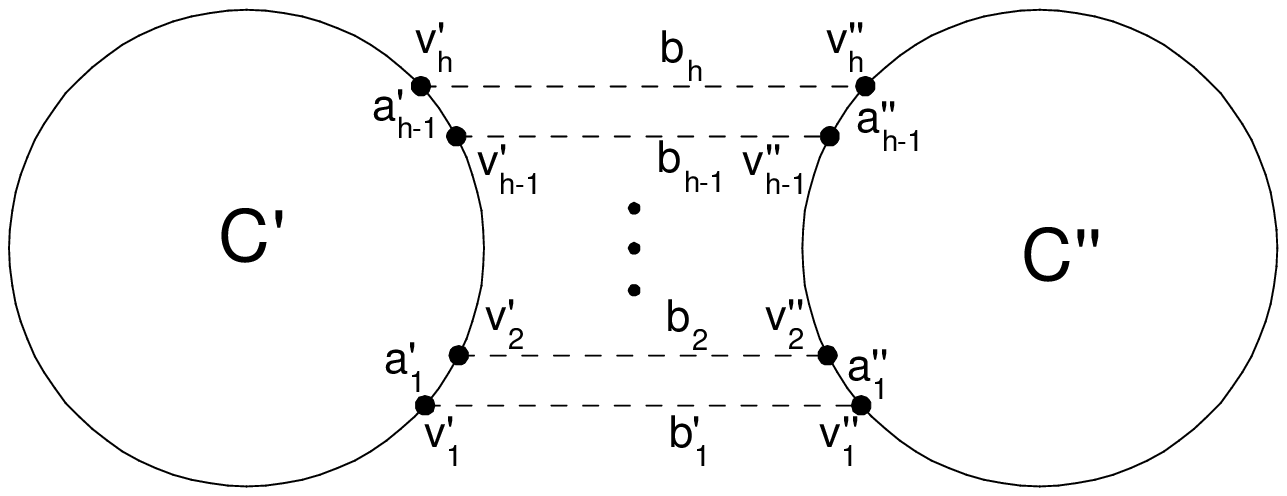}
 \end{center}
 \caption{A $\pqr$-block.}
 \label{Fig. 1}
\end{figure}

The vertices $v_1',v_1'',v_h',v_h''$ are said to be the {\it
corners\/} of the block and the two $\{p,q\}$-residues of $\O$ are
called the {\it sides\/} of the block. We shall often sketch the
block $\O$ as in Figure 2.

If $C'$ and $C''$ are oriented, then a block $\O$ of length $h>1$
is said to be {\it coherent\/} with these  orientations if,
denoted by $v'$ and $w'$ the two corners of $\O$ belonging to $C'$
so that the orientation induced on the side goes from $v'$ to
$w'$, then the orientation induced on the other side goes from
$w''$ to $v''$, where $w''$ and $v''$ are the vertices
$r$-adjacent to $w'$ and $v'$ respectively. In this case the
vertices $v'$ and $w''$ are said to be the {\it key-vertices\/} of
the coherent block (see Figure 2). Each block of length $h=1$ is
considered coherent and both its vertices are key-vertices.

\bigskip

\begin{figure}[bht]
 \begin{center}
 \includegraphics*[totalheight=3.5cm]{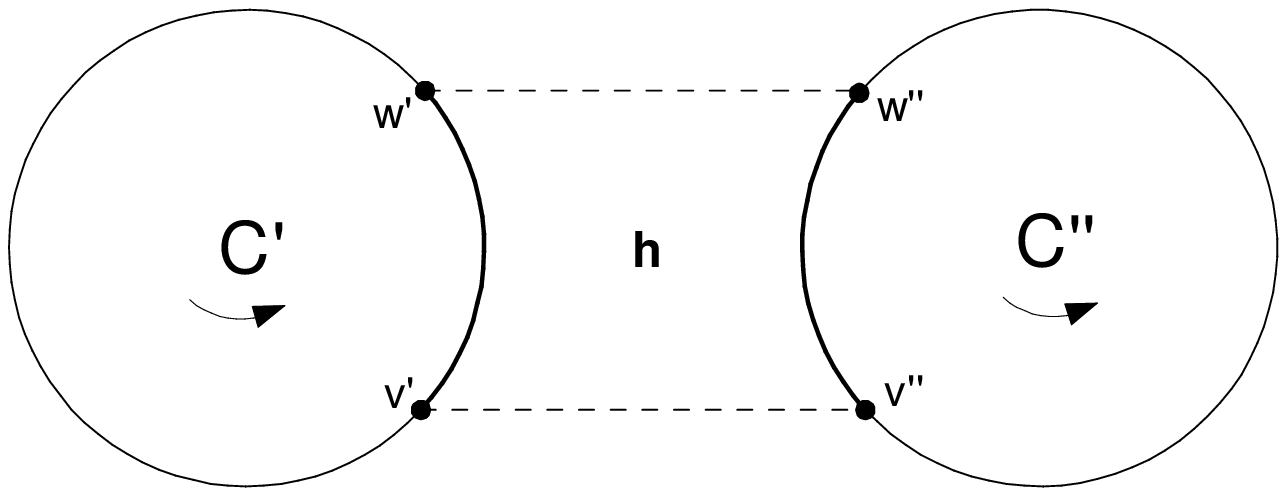}
 \end{center}
 \caption{A coherent block.}
 \label{Fig. 2}

\end{figure}

Suppose now that $\G$ is a $3$-gem, and let $C',C''$ be
$\{p,q\}$-residues of $\G$ belonging to different components of
$\G_{\h r}$. Let $\O$ be a maximal $\pqr$-block, connecting $C'$
with $C''$ and suppose that $\O$ has length $h$. Denote by
$\G(\O)$ the $3$-gem obtained from $\G$ in the following
way:\begin{itemize}
\item delete all the vertices and the edges of $\O$;
\item weld the ``hanging'' edges of the same colour which in $\G$ have $r$-adjacent endpoints belonging to $\O$ (see Figure 3).\end{itemize}

\noindent The graph $\G(\O)$ is said to be {\it obtained by cancelling $\O$ in
$\G$\/}.

\bigskip

\begin{figure}[bht]
 \begin{center}
 \includegraphics*[totalheight=11cm]{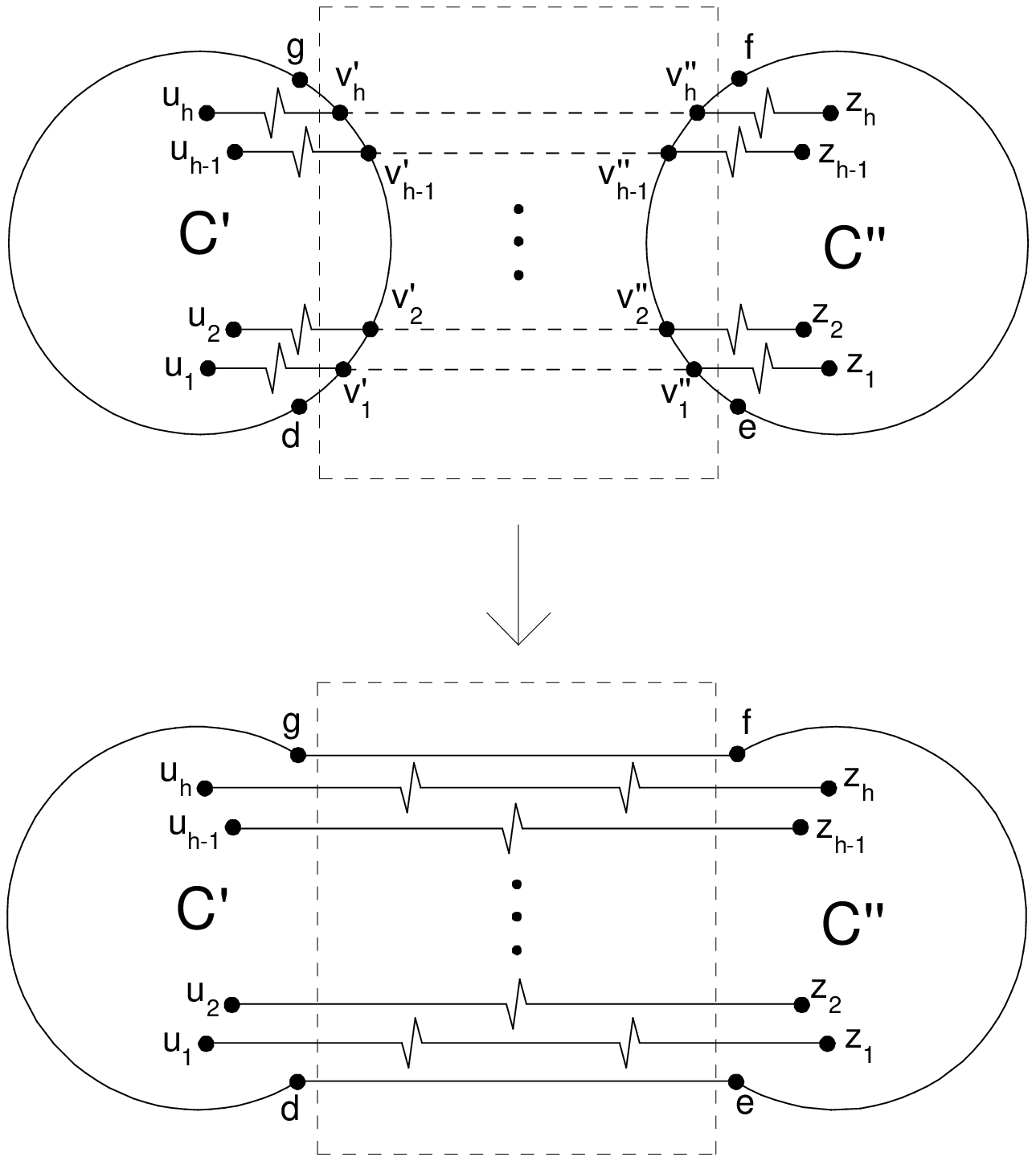}
 \end{center}
 \caption{From $\G$ to $\G(\O)$.}

 \label{Fig. 3}

\end{figure}

\begin{lemma}\label{Lemma 1} The graphs $\G$ and $\G(\O)$ represent the same
$3$-manifold. \end{lemma}

\begin{proof}
Assume for the block $\O$ the notations given in Figure 1. Since
$v_1'$ and $v_1''$ belong to different components of $\G_{\h r}$,
the $r$-edge $\b_1$, together with its endpoints $v_1',v_1''$, is
a dipole of type 1 in $\G$. The cancellation of this dipole
produces a dipole of type 2 involving the $r$-coloured edge
$\b_2$. The sequence of cancellations of this dipole and of the
resulting dipoles of type 2 successively involving
$\b_3,\ldots,\b_h$ leads to $\G(\O)$.
\end{proof}
\medskip

\begin{remark} \label{Remark 1} - It is important to note that $\G$ and
$\G(\O)$ have the same number of $\{r,s\}$-residues, for $s\ne
p,q,r$.
\end{remark}

Notice that $\G(\O)$ is obtained from $\G$ by means of a
``polyhedral gluing'' in the sense of Definition 8 of \cite{FG}. For
this reason, we say that the $\pqr$-block $\O$ is a {\it gluing
subgraph of $\G$} (connecting $C'$ with $C''$ by colour $r$).

After this operation, the $\{p,q\}$-residues $C'$ and $C''$ give
rise to a unique $\{p,q\}$-residue $C$ in $\G(\O)$. Moreover, if
$C'$ and $C''$ are oriented and the block $\O$ is coherent with
these orientations, $C$ inherits an orientation in a natural way.

\section {From $6$-tuples to $3$-manifolds of genus two}

We recall now the possibility of representing all $3$-manifolds of
genus $g\le 2$ via crystallizations defined by $6$-tuples of
integers satisfying simple conditions \cite{CG}.

Let $\wt{\FF}$ be the set of the $6$-tuples
$$f=(h_0,h_1,h_2;q_0,q_1,q_2)$$
of integers satisfying the following conditions:

\begin{itemize}
\item[(I)] $h_i>0$, for each $i\in\Z_3$;
\item[(II)] all $h_i$'s have the same parity;
\item[(III)] $0\le q_i<h_{i-1}+h_i=2l_i$, for each $i\in\Z_3$;
\item[(IV)] all $q_i$'s have the same parity.\end{itemize}
\medskip
\noindent {\it NOTATION\/} - From now on, the operations on the $q_i$ components
will be considered $\text{mod } 2l_i$ and, according to (III), $q_i$ is always the
least non-negative integer of the class.

\medskip

Let $\wt{\GG}=\{\G(f)\mid f\in\wt{\FF}\}$ be the class of
$4$-coloured graphs $\G(f)$ whose vertices are the elements of the
set $$V(f)=\bigcup_{i\in\Z_3}\{i\}\times\Z_{2l_i},$$ and whose
coloured edges are defined by means of the following four
fixed-point-free involutions on $V(f)$:\footnote{Here and in the
following the arithmetic on $V(f)$ is $\text{mod } 3$ in the first
coordinate and $\text{mod } 2l_k$ in the second coordinate of each
vertex $(h,k)\in V(f)$.} \begin{equation}\label{3.1}\i_0(i,j)=(i,j+(-1)^j),\end{equation}
\begin{equation}\label{3.2}\i_1(i,j)=(i,j-(-1)^j),\end{equation}
\begin{equation}\label{3.3}\i_2(i,j)=\begin{cases} (i+1,-j-1) & \text{ if } j=0,\ldots,h_i-1\\
  (i-1,2l_i-j-1)&\text{if} j=h_i,\ldots,2l_i-1\end{cases},\end{equation}
\begin{equation}\label{3.4}\i_3(i,j)=\t\i_2\t^{-1},\end{equation}
 where $\t :V(f)\to V(f)$
is the bijection defined by $$\t(i,j)=(i,j+q_i).$$

\smallskip

To complete the $4$-coloured graph $\G(f)$, join the vertex $v$
with the vertex $w$ by a $c$-coloured edge ($c\in\D_3$) if and
only if $w=\i_c(v)$. Observe that, by (\ref{3.4}), there is a $2$-edge
joining $v_1$ with $v_2$ if and only if there is a $3$-edge
joining $\t(v_1)$ and $\t(v_2)$.

The graph $\G(f)$ contains three $\{0,1\}$-residues $C_i$ of
length $2l_i$, whose vertices are the elements of $V(f)$ having
$i$ as first coordinate. The natural cyclic ordering on
$\Z_{2l_i}$ induces an orientation on each $\{0,1\}$-residue $C_i$
and the bijection $\t$ acts on each $C_i$ as a rotation of
amplitude $q_i$ according to this fixed orientation. Moreover, for
each $i$, there exist a unique maximal $\xd$-block $B_i$ and a
unique maximal $\xt$-block $B_i'$, both of length $h_i$,
connecting $C_i$ with $C_{i+1}$. All these blocks are coherent
with the orientations of the $\xu$-residues.

The map $\t$ is an automorphism of $\G(f)$ exchanging colour $2$
with colour $3$ and, only in case of $q_i$ odd, exchanging colour
$0$ with colour $1$; it is easy to see that $\t$ sends each block
$B_i$ in $B_i'$. Finally, note that $\G(f)$ is bipartite because
of condition (IV).

We shall represent $\G(f)$ by means of a planar embedding of its
residues $\G(f)_{\h 3}$ and $\G(f)_{\h 2}$ (see Figure 4). The
whole graph $\G(f)$ arises by gluing $\G(f)_{\h 3}$ and $\G(f)_{\h
2}$ in the three cycles $C_0,C_1,C_2$, which are the
$\{0,1\}$-residues of $\G(f)$. The graph $\G(f)$ admits a $2$-cell
embedding, which is regular in the sense of \cite{G1}, into an
orientable surface of genus two. This embedding can be obtained in
a standard way using the construction described in \cite{G1}.

\bigskip

\begin{figure}[bht]
 \begin{center}
 \includegraphics*[totalheight=16cm]{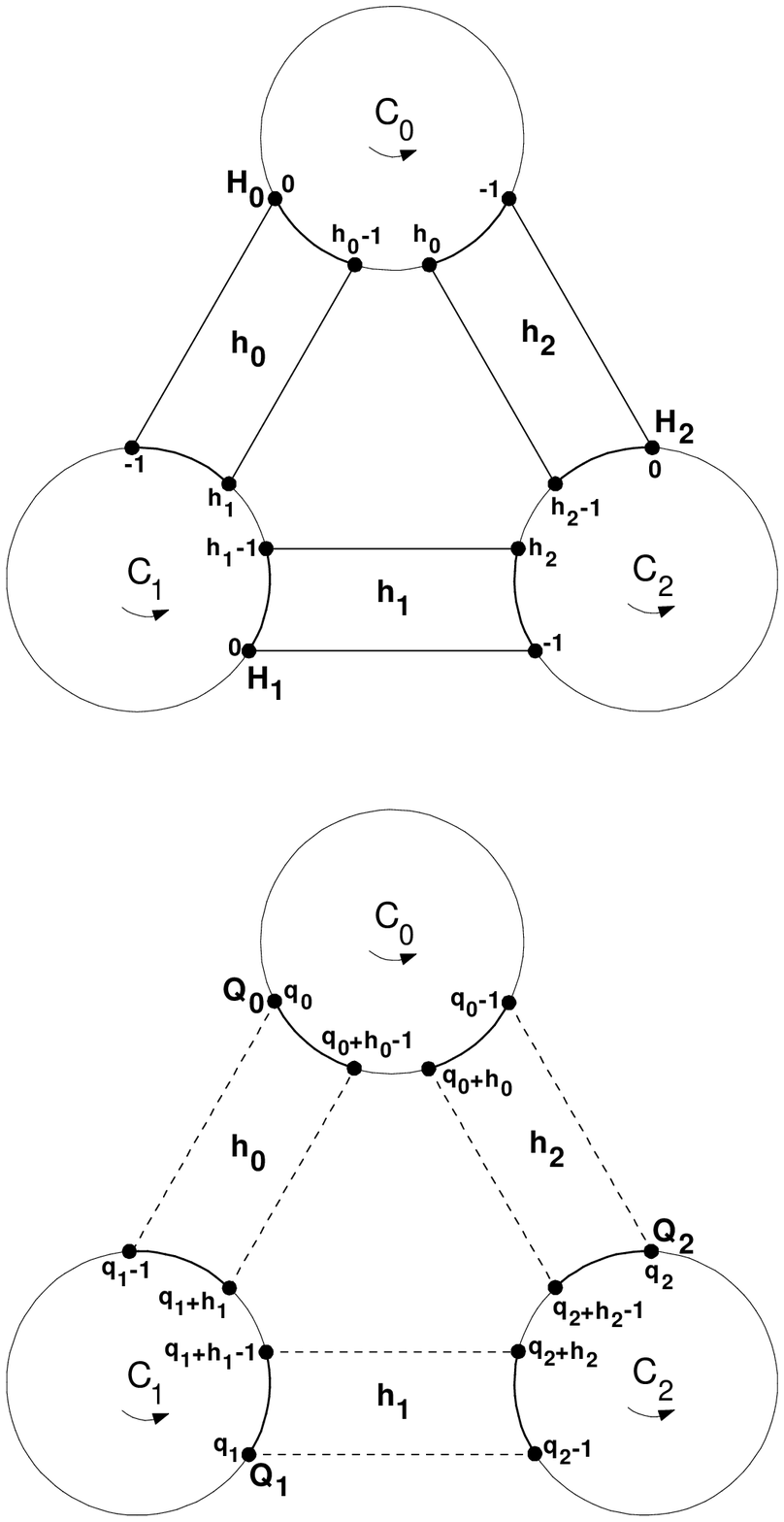}
 \end{center}
 \caption{The 3-residues $\G(f)_{\h 3}$ and $\G(f)_{\h 2}$.}

 \label{Fig. 4}

\end{figure}

We point out that in Figures 4--10 the ``thin'' arcs are edges of
the graph (i.e. there are no vertices in their interior).

\begin{remark} \label{Remark 2} - Suppose now that $\G$ is a bipartite
$4$-coloured graph such that:
\begin{itemize}\item $\G_{\xu}$ consists of
three $\xu$-residues $C_i$ of length $2l_i=h_{i-1}+h_i$, with a
given orientation;
\item for each $i$, there exist a unique maximal $\xd$-block $B_i$
and a unique maximal $\xt$-block, connecting $C_i$ with $C_{i+1}$,
both of length $h_i>0$ and both coherent with the orientations of $C_i$
and $C_{i+1}$.\end{itemize}

\noindent Let $H_i$ (resp. $Q_i$) be the key-vertex of $B_i$
(resp. $B_i'$) belonging to $C_i$. Then $\G$ is the graph
$\G(h_0,h_1,h_2;q_0,q_1,q_2)$, where $q_i$ is the distance from
$H_i$ to $Q_i$ according to the orientation of $C_i$. Observe that
all $q_i$'s have the same parity since $\G$ is bipartite.
\end{remark}

Denote by $\FF$ the subset of $\wt{\FF}$ consisting of the
$6$-tuples $f$ such that:

\begin{itemize}
\item[(V)] $h_i+q_i$ is odd, for each $i\in\Z_3$;
\item[(VI)] $\G(f)$ contains exactly three $\{2,3\}$-residues.
\end{itemize}

Then $\G(f)\in\wt{\GG}$ is a crystallization of a $3$-manifold of
genus $g\le 2$ if and only if $f\in\FF$ \cite{CG}.
We shall call {\it admissible\/} each $6$-tuple belonging to $\FF$.
Observe that, as a consequence of (II) and (V), $q_i\ne h_j$ for each $i,j\in\Z_3$.

Set $\GG=\{\G(f)\mid f\in\FF\}$; the crystallizations of $\GG$ are
{\it $2$-symmetric\/} in the sense of \cite{CG} and hence they represent
$2$-fold branched coverings of $S^3$. Moreover, if $M$ is a
$3$-manifold of genus $g\le 2$, then there exists an $f\in\FF$
such that $\G(f)\in\GG$ is a crystallization of $M$ \cite{CG}. Thus, the
set of all admissible $6$-tuples gives a {\it complete catalogue
of all $3$-manifolds of genus $g\le 2$} (see \cite{Ca}).

The open problem of classifying $3$-manifolds of genus two can be translated
into the following question: when do two admissible $6$-tuples represent
the same manifold?

In this direction, it is important to find ``elementary
transformations'' on admissible $6$-tuples, which change the
associated graph but do not change the represented manifold. The
present paper describes an elementary transformation of this type,
which is called {\it $2$-symmetric\/}, since it can be obtained by
considering standard sequences of dipole moves which change a
given $2$-symmetric crystallization $\G(f)\in\GG$ to another
$2$-symmetric crystallization $\G'=\G(f')\in\GG$. The induced
$2$-symmetric transformation changes the admissible $6$-tuple $f$
into the admissible $6$-tuple $f'$ representing the same manifold.

We claim that for particular values of the parameters of $f$, the
graph $\G(f)$ represents a $3$-manifold of genus $0$ or $1$.

\begin{lemma} \label{Lemma 2}Let $\hq$ be an admissible $6$-tuple and let
$(i,j,k)$ be any permutation of $\Z_3$.
\begin{itemize}
\item[a) ] If $q_i=q_j=0$ then $\G\hq$ represents the lens space $L(l_k,q_k/2)$.
\item[b) ] If $q_0=q_1=q_2=0$ then $\G\hq$ represents $\mathbf S^3$.\end{itemize}
\end{lemma}

\begin {proof} a) By deleting all dipoles involving the $2$- and
$3$-edges connecting $C_i$ with $C_j$ we obtain the ``normal''
crystallization of the lens space $L(l_k,q_k/2)$ (see \cite{DG}). b)
$L(l_k,0)\cong\mathbf S^3$.
\end {proof}

\section { The 2-symmetric transformation}

Let $\G(f)\in\GG$ be the crystallization of a $3$-manifold $M$
defined by the $6$-tuple $f=\hq\in\FF$. Moreover, assume the
notation of the previous section and suppose the cycles $C_i$
oriented according to the natural cyclic ordering on $\Z_{2l_i}$
(see Figure 4).

Delete the following edges from $\G(f)$:
\begin{itemize}
\item all the $h_1$ $2$-edges connecting $C_1$ with $C_2$;
\item the $1$-edge of $C_0$ connecting $(0,-1)$ with $(0,0)$;
\item the edge $\a$ of $C_0$ connecting $(0,h_0-1)$ with $(0,h_0)$.\footnote{Note that $\a$ has colour 0 (resp. 1) if $h_0$ is odd (resp. even).}\end{itemize}

\noindent Denote by $\wt{\G}(f)$ the resulting $4$-coloured graph with boundary.

Let now $\G(h_1)$ be the $4$-coloured graph with boundary shown in Figure 5.

\bigskip

\begin{figure}[bht]
 \begin{center}
 \includegraphics*[totalheight=5cm]{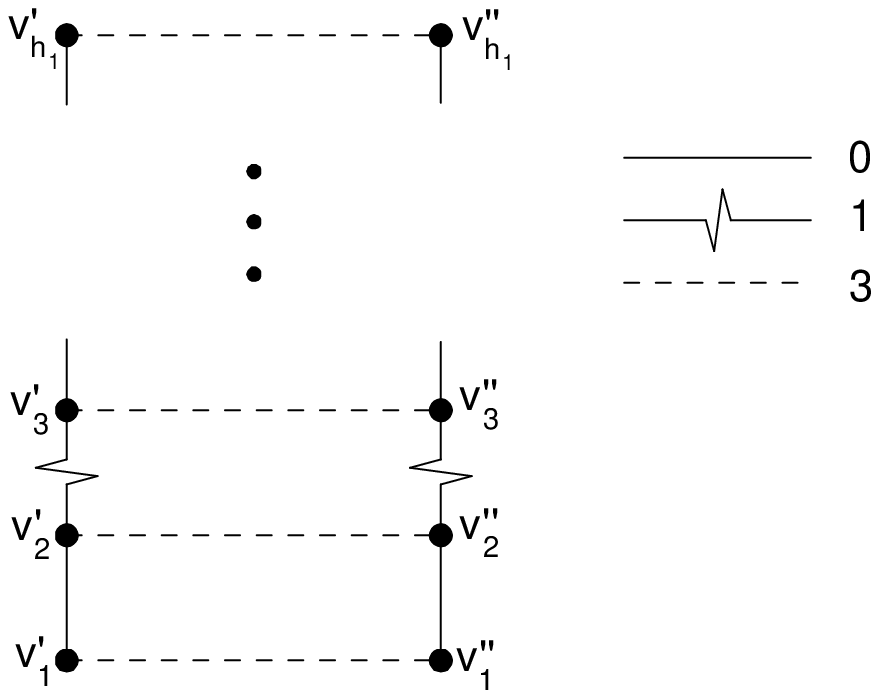}
 \end{center}
 \caption{$\G(h_1)$.}

 \label{Fig. 5}

\end{figure}

\begin{remark} \label{Remark 3}  - All vertices of $\G(h_1)$ are boundary
vertices with respect to colour $2$. Moreover, $v_1'$ and $v_1''$
are boundary vertices with respect to colour $1$, and $v_{h_1}'$
and $v_{h_1}''$ are boundary vertices with respect to colour $0$
(resp. $1$) if the $h_i$'s are odd (resp. even), i.e. if $\a$ has
colour $0$ (resp. $1$).
\end{remark}

Now connect:
\begin{itemize}
\item the vertex $(0,0)$ (resp. $(0,-1)$) of $\wt{\G}(f)$ with the
vertex $v_1'$ (resp. $v_1''$) of $\G(h_1)$ by a $1$-edge;
\item the vertex $(0,h_0-1)$ (resp. $(0,h_0)$) of $\wt{\G}(f)$ with the
vertex $v_{h_1}'$ (resp. $v_{h_1}''$) of $\G(h_1)$ by a
$\g(\a)$-edge (recall the previous remark);
\item the vertices $(1,0),\ldots,(1,h_1-1)$ of $\wt{\G}(f)$ respectively
with the vertices $v_1',\ldots,v_{h_1}'$ of $\G(h_1)$ by $2$-coloured
edges;
\item the vertices $(2,-1),(2,-2),\ldots,(2,-h_1=h_2)$ of $\wt{\G}(f)$
respectively with the vertices $v_1'',v_2''\ldots,v_{h_1}''$ of $\G(h_1)$
by $2$-coloured edges.\end{itemize}

\noindent Denote by $G(f)$ the resulting $4$-coloured graph (without boundary).

\bigskip

\begin{figure}[bht]
 \begin{center}
 \includegraphics*[totalheight=12cm]{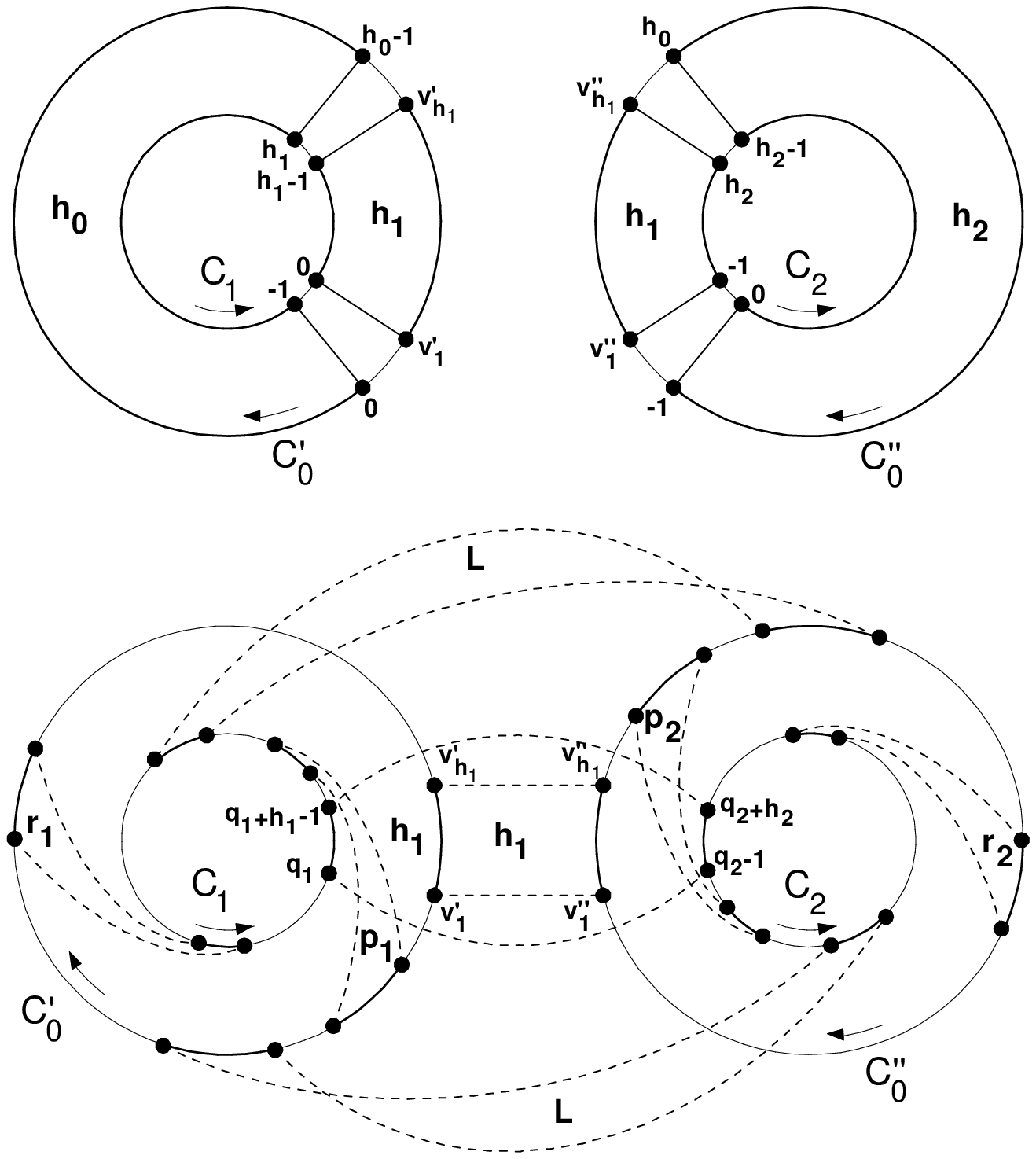}
 \end{center}
 \caption{The 3-residues $G(f)_{\h 3}$ and $G(f)_{\h 2}$.}

 \label{Fig. 6}

\end{figure}

Note that $\G(h_1)$ is the subgraph of $G(f)_{\h 2}$ induced by
the set of vertices $\{v_1',\ldots,v_{h_1}',v_1'',\ldots,v_{h_1}''\}$.

The $\{0,1\}$-residue $C_0$ of $\G(f)$ splits in $G(f)$ into two different
components $C_0'$ and $C_0''$, where:
\begin{itemize}
\item the sequence of the vertices $(0,0),\ldots,(0,h_0-1)$ of $\G(f)$
followed by the sequence of the vertices $v_{h_1}',v_{h_1-1}',\ldots,v_1'$
of $\G(h_1)$ gives all consecutive vertices of $C_0'$;
\item the sequence of the vertices $(0,h_0),(0,h_0-1)\ldots,(0,-1)$ of $\G(f)$
followed by the sequence of the vertices $v_1'',\ldots,v_{h_1}''$ of $\G(h_1)$
gives all consecutive vertices of $C_0''$.\end{itemize}

The graph $G(f)_{\h 3}$ has two components $\O'$ and $\O''$; $\O'$
(resp. $\O''$) has two $\xu$-residues $C_0',C_1$ (resp.
$C_0'',C_2$) of length $h_0+h_1$ (resp. $h_1+h_2$), connected by
``parallel'' $2$-edges. Hence, the $\xt$-block $\G(h_1)$ is a
gluing subgraph of $G(f)$ (connecting $C_0'$ with $C_0''$ by
colour $3$) of length $h_1$. Moreover, the graph obtained by
cancelling $\G(h_1)$ in $G(f)$ is $\G(f)$. This proves that $G(f)$
represents the $3$-manifold $M$.

We are now going to show that, by choosing another suitable gluing
subgraph of $G(f)$, we can obtain a new crystallization $\G(f')\in\GG$
of the $3$-manifold $M$, depending on a different $6$-tuple $f'\in\FF$.
To achieve this goal, relabel the vertices of $C_0',C_0''$
and $C_2$ of $G(f)$ in the following way:
\begin{itemize}
\item label the vertices of $C_0'$ by $(0',j)$, $j\in\Z_{2l_1-1}$, so that
in the increasing sequence $(0',0),\ldots,(0',2l_1-1)$ the vertices are
consecutive and so that the vertices $v_i'$ of $C_0'$ are labelled by
$(0', q_1+h_1-i)$, for each $i=1,\ldots,h_1$;
\item label the vertices of $C_0''$ by $(0'',j)$, $j\in\Z_{2l_2-1}$, so that
in the increasing sequence $(0'',0),\ldots,(0'',2l_2-1)$ the vertices are
consecutive and so that the vertices $v_i''$ of $C_0''$ are labelled by
$(0'', q_2+i-1)$, for each $i=1,\ldots,h_1$;
\item relabel the vertices of $C_2$ so that the second component of $(2,j)$ becomes $(2,j-h_2)$,
for each $j\in\Z_{2l_2-1}$.\end{itemize}

\bigskip

\begin{figure}[bht]
 \begin{center}
 \includegraphics*[totalheight=12cm]{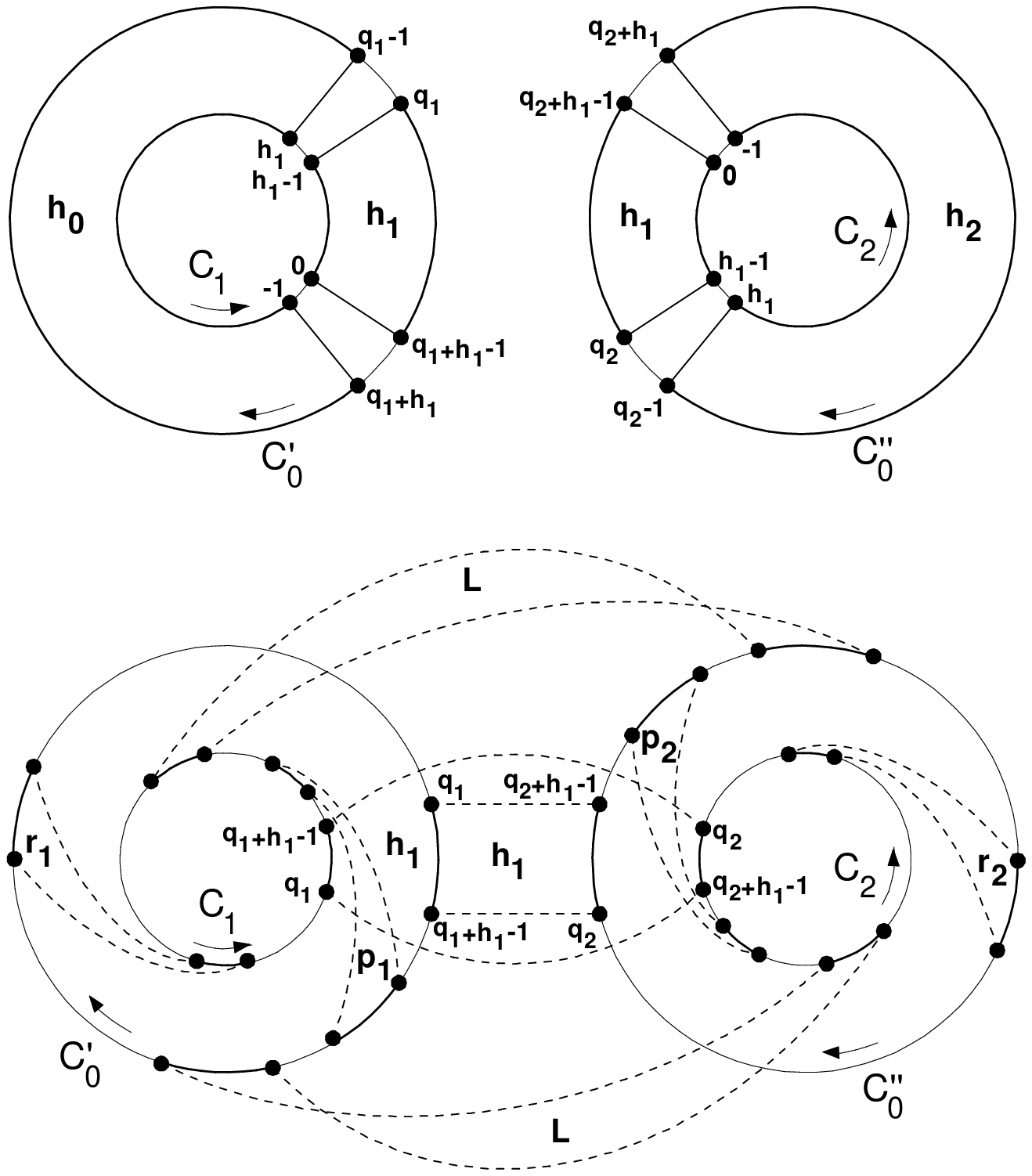}
 \end{center}
 \caption{The 3-residues $G(f)_{\h 3}$ and $G(f)_{\h 2}$.}

 \label{Fig. 7}

\end{figure}

Assume on $C_0'$ and $C_0''$ the orientations induced by the cyclic ordering of their vertex labellings:
the gluing subgraph $\G(h_1)$ is coherent with these orientations and its cancellation restores the original
orientation on $C_0$.

\begin{remark} \label{Remark 4}  - The subgraph $\G'(h_1)$ of $G(f)_{\h 2}$
induced by the set of vertices
$\{(1,q_1),\ldots,(1,q_1+h_1-1),(2,q_2),\ldots,(2,q_2+h_1-1)\}$ is
a gluing subgraph of $G(f)$, connecting $C_1$ with $C_2$ by colour
$3$. The $4$-coloured graph obtained by cancelling $\G'(h_1)$ in
$G(f)$ is c.p.-isomorphic to $\G(f)$. This follows immediately
since the involution on $V(G(f))$ exchanging $(1,i)$ with
$(0',i)$, for each $i=0,1,\ldots,h_0+h_1$, and $(2,j)$ with
$(0'',j)$, for each $j=0,1,\ldots,h_1+h_2$ is a c.p.-automorphism
of $G(f)$ sending $\G(h_1)$ to $\G'(h_1)$.
\end{remark}

Let now $\TE$ (resp. $\TE'$) denote the unique gluing subgraph
of $G(f)$ connecting $C_1$ with $C_0''$ (resp. $C_0'$ and $C_2$) by colour
$3$. As can be easily checked, $\TE$ and $\TE'$ are nonvoid if and only
if $q_0\ne 0$.
The involutory c.p.-automorphism defined in previous remark
sends $\TE$ to $\TE'$ and therefore the $4$-coloured graphs respectively
obtained by deleting $\TE$ and $\TE'$ in $G(f)$ are c.p.-isomorphic.
From now on, we focus our attention on the $4$-coloured graph $\G'$
obtained by cancelling $\TE$ in $G(f)$. In fact, this is the unique
graph obtained by cancelling gluing subgraphs in $G(f)$ connecting
$\{0,1\}$-residues by colour $3$, which is, in general, different from
$\G(f)$, up to c.p.-isomorphisms. It is straightforward that $\G'$
still represents the $3$-manifold $M$. Moreover, the following result
holds:

\begin{theorem}\label{Theorem 3}
Let $f=\hq$ be an admissible $6$-tuple such that $q_0\ne 0$
and let $f'=\hqp$ be the $6$-tuple defined by the following rules:
\begin{equation}\label{4.1}\begin{cases}
h_0'=h_0+h_1-q_0\\h_1'=q_0\\h_2'=h_2+h_1-q_0\end{cases}\quad\begin{cases}
q_0'=h_0+h_1+h_2-2q_0\\q_1'=q_0+q_1+h_1\\q_2'=q_0+q_2+h_1\end{cases},\quad\text{if}\,\,\,\,0<q_0<h_0,h_2;\end{equation}
\begin{equation}\label{4.2}\begin{cases}
h_0'=q_0+h_1-h_2\\h_1'=h_0+h_2-q_0\\h_2'=q_0+h_1-h_0\end{cases}\quad\begin{cases}
q_0'=h_1\\q_1'=q_0+q_1-h_2\\q_2'=q_0+q_2-h_0\end{cases},\quad\text{if}\,\,\,\,q_0>h_0,h_2;\end{equation}
\begin{equation}\label{4.3}\begin{cases}
h_0'=h_1\\h_1'=h_0\\h_2'=h_1+h_2-h_0\end{cases}\quad\begin{cases}
q_0'=h_1+h_2-q_0\\q_1'=q_1\\q_2'=2q_0+q_2+h_1-h_0\end{cases},\quad\text{if}\,\,\,\,h_0<q_0<h_2;\end{equation}
\begin{equation}\label{4.4}\begin{cases}
h_0'=h_1+h_0-h_2\\h_1'=h_2\\h_2'=h_1\end{cases}\quad\begin{cases}
q_0'=h_1+h_0-q_0\\q_1'=2q_0+q_1+h_1-h_2\\q_2'=q_2\end{cases},\quad\text{if}\,\,\,\,h_2<q_0<h_0.\end{equation}
Then $f'$ is an admissible $6$-tuple and the $4$-coloured graphs
$\G(f)$ and $\G(f')$ represent the same manifold.
\end{theorem}

\begin {proof} With the previous assumptions and notation, it suffices to show that the $4$-coloured graph
$\G'$ obtained by cancelling $\TE$ in $G(f)$ is c.p.-isomorphic to $\G(f')$.

First of all, exchange the names of the two cycles $C_0'$ and $C_1$, together with the first components
in the labelling of their vertices. After this relabelling, $\TE$ becomes the unique gluing subgraph of
$G(f)$ connecting $C_0'$ with $C_0''$ by colour $3$; denote by $L$ the length of $\TE$.

Figure 8 sketches, with the usual conventions, the graphs
$G(f)_{\h 3}$ and $G(f)_{\h 2}$; here we also point out the
labelling of some ``strategic'' vertices of $G(f)$. The
computation of the integers $L,p_1,p_2,r_1,r_2$, depending on the
components of $f$, is described in the following table:

\smallskip

\vbox{
\tabskip=1pc
%\moveright 2 in
\vbox{\offinterlineskip
\halign{\strut#\hfil&\vrule\quad\hfil#\hfil&\hfil#\hfil&\hfil#\hfil&\hfil#\hfil&\hfil#\hfil\quad\vrule \cr
&$L$&$p_1$&$p_2$&$r_1$&$r_2$\cr
\noalign{\hrule}
if $\,\,0<q_0<h_0,h_2$&$q_0$&$0$&$0$&$h_0-q_0$&$h_2-q_0$\cr
if $\,\,q_0>h_0,h_2$&$h_0+h_2-q_0$&$q_0-h_2$&$q_0-h_0$&$0$&$0$\cr
if $\,\,h_0<q_0<h_2$&$h_0$&$0$&$q_0-h_0$&$0$&$h_2-q_0$\cr
if $\,\,h_2<q_0<h_0$&$h_2$&$q_0-h_2$&$0$&$h_0-q_0$&$0$\cr
\noalign{\hrule}
}
}

\smallskip

\centerline{Table 1}
}

\smallskip

Note that:\begin{itemize}
\item[(A)] $p_1\ne0$ if and only if $r_2=0$;
\item[(B)] $p_2\ne0$ if and only if $r_1=0$ \end{itemize}

Now, we are going to look into the shape of $\G_{\h 3}'$ and
$\G_{\h 2}'$. It is clear that, by cancelling $\TE$ in $G(f)$, the
two $\{0,1\}$-cycles $C_0'$ and $C_0''$ of $G(f)$ give rise to a
unique $\{0,1\}$-cycle $C_0$ of $\G'$; moreover, since the length
of $C_0'$ (resp. $C_0''$) is $h_0+h_1$ (resp. $h_1+h_2$), the
length of $C_0$ is $h_0+2h_1+h_2-2L$. Since the gluing subgraph
$\TE$ is coherent with the orientations on $C_0'$ and $C_0''$, the
cycle $C_0$ inherits an orientation in a natural way. On the other
hand, the length of $C_1$ (resp. $C_2$) in $\G'$ is still
$h_0+h_1$ (resp. $h_1+h_2$). Hence, the $3$-coloured graphs
$\G_{\h 3}'$ and $\G_{\h 2}'$ can be sketched as in Figure 9 and
10 respectively.

By properties (A) and (B), it is easy to check that, in all four
cases of Table 1, the graph $\G_{\h 2}'$ is planar and has the
shape of Figure 9, where the numbers inside the strips can be
computed by Table 1.

The graph $\G'$ satisfies the assumptions of Remark \ref{Remark 2}
and hence $\G'$ is c.p.-isomorphic to $\G(f')\in\wt{\GG}_2$, where
$f'=(h_0+h_1-L,L,h_1+h_2-L;q_0',q_1',q_2')$.

We are now going to compute $q_1',q_2'$ and $q_0'$. If $H_i$
(resp. $Q_i$) denotes the key-vertex of $B_i$ (resp. $B_i'$)
belonging to $C_i$, $i=0,1,2$, then $q_i'$ is the distance from
$H_i$ to $Q_i$ according to the orientation of $C_i$. Now, $H_1$
(resp. $H_2$) is the vertex of $G(f)$ which is $2$-adjacent with
$(0',q_1+h_1+p_1+L-1)$ (resp. with $(0'',q_2+h_1+p_2-1)$) in
$G(f)$. On the other hand, the vertex which is $2$-adjacent with
$(1,0)$ (resp. with $(2,0)$) in $G(f)$ is $(0',q_1+h_1-1)$ (resp.
$(0'',q_2+h_1-1)$); hence, the distance from $H_1$ to  $(1,0)$
(resp. from $H_2$ to $(2,0)$), according to the orientation of
$C_1$ (resp. $C_2$), equals the distance from $(0',q_1+h_1-1)$ to
$(0',q_1+h_1+p_1+L-1)$ (resp. from $(0'',q_2+h_1-1)$ to
$(0'',q_2+h_1+p_2-1)$), according to the orientation of $C_0'$
(resp. $C_0''$). Since $Q_1$ (resp. $Q_2$) is the vertex
$(1,q_1+h_1+p_1)$ (resp. $(2,q_2+h_1+p_2+L)$) of $G(f)$, we
obtain: $$q_1'=(p_1+L)+(q_1+h_1+p_1)=2p_1+q_1+h_1+L,$$
$$q_2'=p_2+(q_2+h_1+p_2+L)=2p_2+q_2+h_1+L.$$ Furthermore, $H_0$ is
the vertex of $C_0$ which is $2$-adjacent in $\G'$ with the vertex
preceding $H_1$ in $C_1$; hence, $H_0$ is the vertex
$(0',q_1+h_1+p_1+L)$ in $C_0$. In the same way, $Q_0$ is the
vertex of $C_0$ which is $3$-adjacent with $(1,q_1+h_1+p_1-1)$ in
$\G'$. Therefore, we have the following two possibilities:
$$Q_0=\begin{cases} (0',q_1+h_1)&\text{ if } p_1\ne 0\\
(0'',q_2)&\text{ if } p_1=0\end{cases}.$$ By recalling that
$p_1\ne 0$ if and only if $r_2=0$, we can conclude that, in both
cases $$q_0'=h_1+r_1+r_2.$$ The graph $\G'$ is c.p.-isomorphic to
$\G(f')\in\wt\GG$, where:
$$f'=(h_0+h_1-L,L,h_1+h_2-L;h_1+r_1+r_2,2p_1+q_1+h_1+L,2p_2+q_2+h_1+L).$$
Substituting in this expression the values of $L,p_1,p_2,r_1$ and
$r_2$ of Table 1 we obtain (\ref{4.1}), (\ref{4.2}), (\ref{4.3})
and (\ref{4.4}); moreover, $f'$ satisfies property (V) and, by
Remark \ref{Remark 1}, property (VI). So, $f'$ is an admissible
$6$-tuple and this completes the proof.
\end {proof}

\bigskip

\begin{figure}[bht]
 \begin{center}
 \includegraphics*[totalheight=12cm]{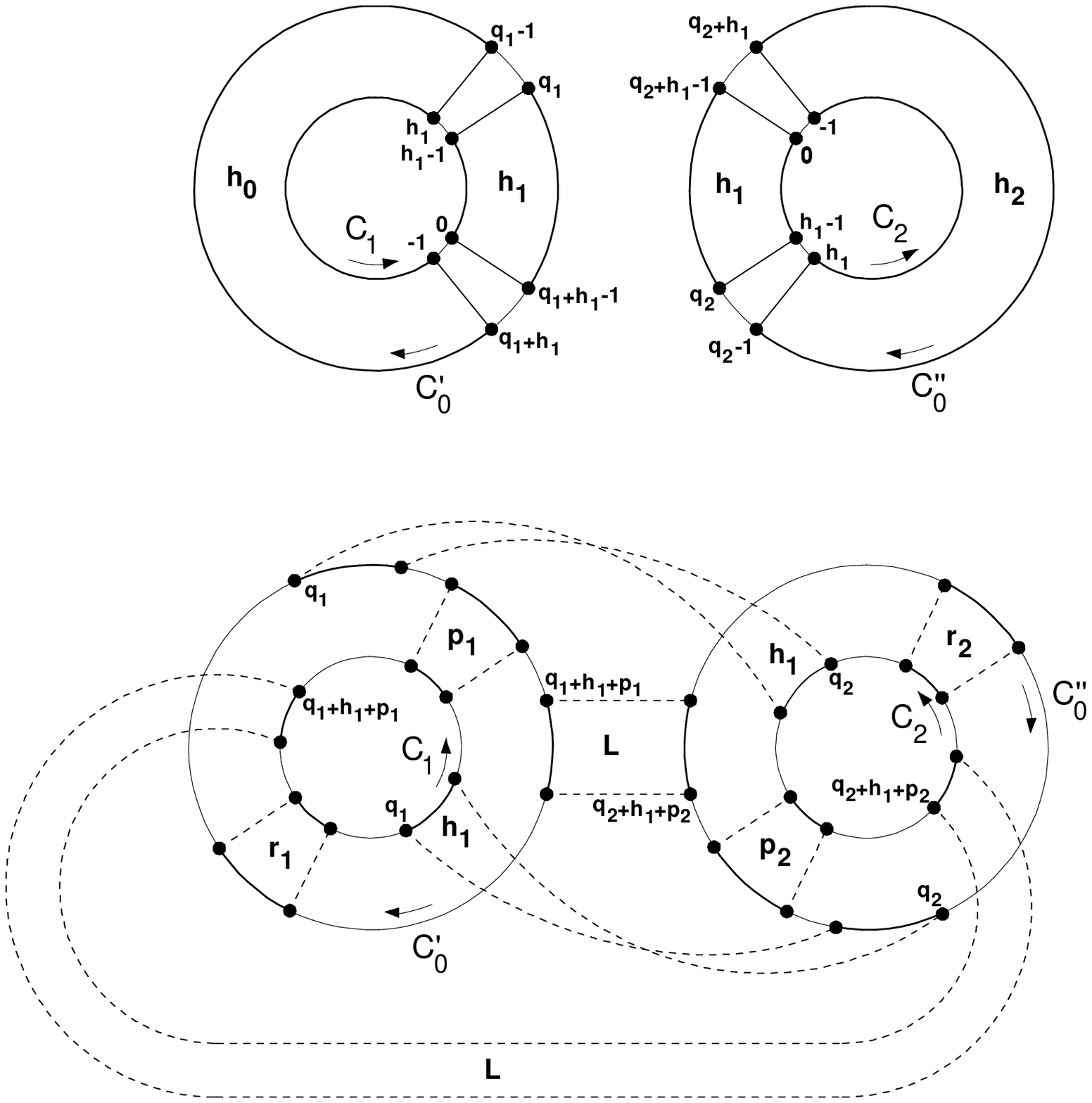}
 \end{center}
 \caption{The 3-residues $G(f)_{\h 3}$ and $G(f)_{\h 2}$.}

 \label{Fig. 8}

\end{figure}

\bigskip

\begin{figure}[bht]
 \begin{center}
 \includegraphics*[totalheight=8cm]{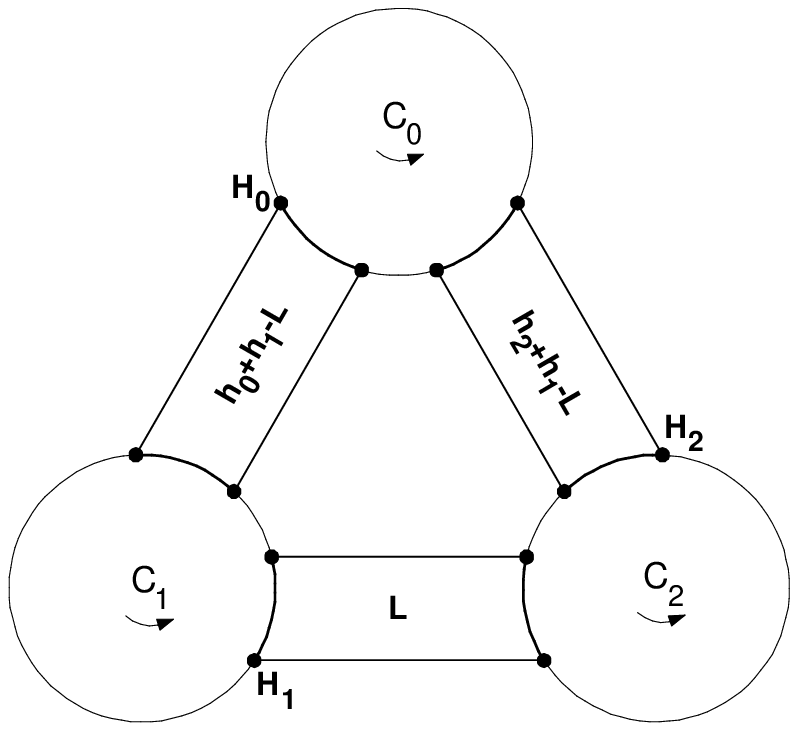}
 \end{center}
 \caption{The 3-residue $\G_{\h 3}'$.}

 \label{Fig. 9}

\end{figure}

\bigskip

\begin{figure}[bht]
 \begin{center}
 \includegraphics*[totalheight=10.5cm]{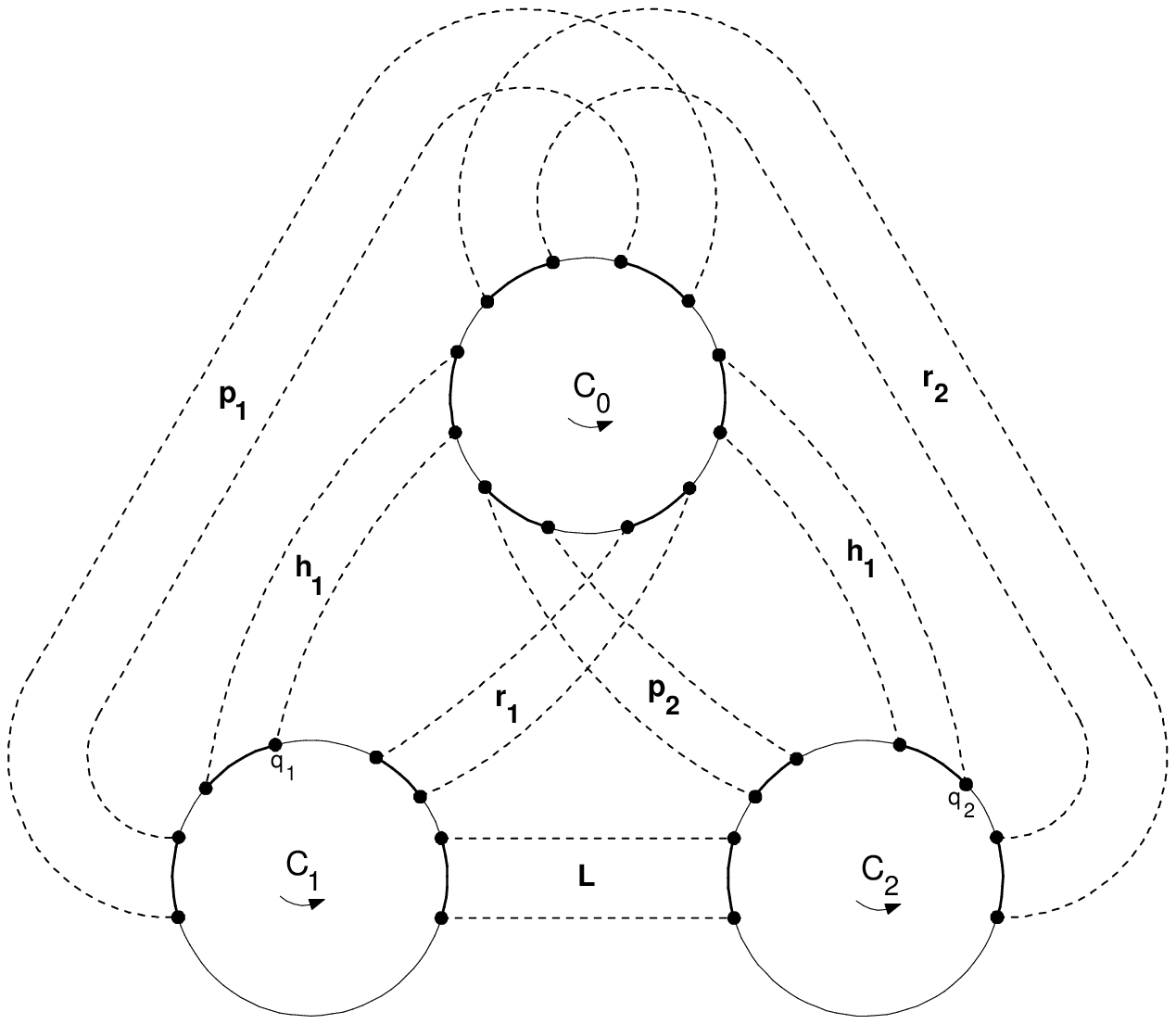}
 \end{center}
 \caption{The 3-residue $\G_{\h 2}'$.}

 \label{Fig. 10}

\end{figure}

With the assumptions of Theorem \ref{Theorem 3}, the
transformation changing $f$ into $f'$ is said to be a {\it
$2$-symmetric transformation\/}.

\section {Equivalence of admissible $6$-tuples}

The reader might suspect that different admissible $6$-tuples can
be associated to c.p.-isomorphic coloured graphs. This is true,
since we can change the order of the three $\{0,1\}$-residues or
their orientations and this choice leads to different $6$-tuples
arising from the same graph.

\begin{lemma}\label{Lemma 4} If $f=\hq$ is an admissible $6$-tuple, then the
$6$-tuples $(h_1,h_2,h_0;q_1,q_2,q_0)$,
$(h_2,h_1,h_0;q_0,q_2,q_1)$, $(h_0,h_1,h_2;-q_0,-q_1,-q_2)$ are
admissible and their associated graphs are c.p.-isomorphic to
$\G(f)$.
\end{lemma}

\begin {proof} See \cite{Ca}, Proposition 16.
\end {proof}

Let $\z_1,\z_2,\z_3:\FF\to\FF$ be the relative maps on the set of all admissible $6$-tuples:
$$\z_1\hq=(h_1,h_2,h_0;q_1,q_2,q_0)$$
$$\z_2\hq=(h_2,h_1,h_0;q_0,q_2,q_1)$$
$$\z_3\hq=(h_0,h_1,h_2;-q_0,-q_1,-q_2)$$

These maps are bijections on $\FF$ such that
$\z_1^3=\z_2^2=\z_3^2=1$. Each of them sends an admissible
$6$-tuple to a (generally different) admissible $6$-tuple
associated to a c.p.-isomorphic graph.

\begin{remark} \label{Remark 5}  - We can interpret the action of $\z_3$ as a
change of orientation of the three $\{0,1\}$-residues $C_0,C_1$
and $C_2$, the action of $\z_1$ as a cyclic permutation $C_0\to
C_1\to C_2\to C_0$ and the action of $\z_2$ as an exchange between
$C_1$ and $C_2$. \end{remark}

Let $\text{Aut}(\FF)$ be the group of all bijections of $\FF$ and
let $K$ be any subgroup of $\text{Aut}(\FF)$; then two admissible
$6$-tuples $f,f'$ will be called {\it $K$-equivalent\/} if there
exists $k\in K$ such that $f=k(f')$. As usual, we call a {\it
$K$-orbit\/} any $K$-equivalence class of admissible $6$-tuples,
i.e. any element of $\FF/K$.

Now, let $H$ and $H'$ be the following subgroups of
$\text{Aut}(\FF)$: $$H=<\z_1,\z_2,\z_3>,\quad H'=<\z_2,\z_3>.$$

\begin{lemma} \label{Lemma 5} The group $H'$ is isomorphic to the Klein four
group $\Z_2\oplus\Z_2$ and the group $H$ is isomorphic to the
dihedral group $D_6$ of all symmetries of a regular hexagon.
Moreover, $H=H'\cup H'\z_1\cup H'\z_1^2$.
\end{lemma}

\begin {proof} The relations $\z_1^3=\z_2^2=\z_3^2=1$,
$\z_2\z_3=\z_3\z_2$,  $\z_1\z_3=\z_3\z_1$ and
$\z_1\z_2=\z_2\z_1^2$ hold. Therefore, we get for $H'$ the
classical presentation of the Klein four group:
$H'\cong<\z_2,\z_3\mid\z_2^2,\z_3^2,[\z_2,\z_3]>\cong\Z_2\oplus\Z_2$.
To obtain the second result, define $r=\z_1^{-1}\z_3,s=\z_3\z_2$
and observe that both $r$ and $s$ commute with $\z_3$. By Tietze
transformations, we have:

\smallskip
%\begin{itemize}

\noindent $\,\,H\cong<\z_1,\z_2,\z_3\mid\z_1^3,\z_2^2,\z_3^2,(\z_1\z_2)^2,[\z_1,\z_3],[\z_2,\z_3]>$,

\noindent $\,\,H\cong<\z_3,r,s\mid(\z_3r^{-1})^3,(\z_3s)^2,\z_3^2,(\z_3r^{-1}\z_3s)^2,\z_3r^{-1}\z_3r,\z_3s\z_3s^{-1}>$,

\noindent $\,\,H\cong<\z_3,r,s\mid \z_3r^{-3},s^2,\z_3^2,(r^{-1}s)^2,\z_3^2,\z_3^2>$,

\noindent $\,\,H\cong<r,s\mid s^2,r^6,(sr)^2>$,

%\end{itemize}
\smallskip

\noindent which is a classical presentation of $D_6$. Finally, the
last sentence holds since $\vert H:H'\vert =\vert <\z_1>\vert=3$
and $H'\cap <\z_1>=\{1\}$.
\end {proof}

Now, let $\FF_H=\FF/H$; then each orbit of $\FF_H$ is composed by
12 (not necessarily distinct) admissible $6$-tuples associated to
c.p.-isomorphic 4-coloured graphs.

The {\it complexity\/} of an admissible $6$-tuple $\hq$ is the
integer $$\v\hq=h_0+h_1+h_2,$$ which is half the cardinality of
$V(\G(f))$. Since $6$-tuples of the same $H$-orbit have the same
complexity, we can translate the notion of complexity to
$H$-orbits in an obviuos way.

To avoid repetitions of c.p.-isomorphic graphs, it is very useful
to select a canonical representative for each $H$-orbit.

\begin{lemma} \label{Lemma 6} If $\o$ is an $H$-orbit, then there exists a
unique $6$-tuple $f=\hq\in\o$ such that the following conditions
hold:
\begin{itemize}
\item[(a) ] $h_0\le h_1\le h_2$;
\item[(b) ] $q_0\le l_0$;
\item[(c) ] if $q_0=0,l_0$ then $q_1\le l_1$;
\item[(d) ] if $q_0=0,l_0$ and $q_1=0,l_1$ then $q_2\le l_2$;
\item[(e) ] if $h_0=h_1$ then $q_0\le q_2$ and $q_2\le -q_0$;
\item[(f) ] if $h_0=h_1$ and $q_2=\pm q_0$ then $q_1\le h_1$;
\item[(g) ] if $h_1=h_2$ then $q_0\le q_1$ and $q_1\le -q_0$;
\item[(h) ] if $h_1=h_2$ and $q_1=\pm q_0$ then $q_2\le h_2$;
\item[(i) ] if $h_0=h_1=h_2$ then $q_1\le q_2$.
\end{itemize}
\end{lemma}

\begin {proof}
By $\z_1$ and $\z_2$ we can permute $h_0,h_1,h_2$ in
all possible way and therefore condition (a) can be achieved.
Conditions (b), (c) and (d) follow by a suitable application of
$\z_3$. Conditions (e),(f),(g),(h) and (i) follow by a combined application
of the three maps. The unicity of such an $f$ is
straightforward.
\end {proof}

The $6$-tuple $f$ of the previous lemma is said to be the {\it
canonical representative\/} of the $H$-orbit $\o$.

\begin{remark}  \label{Remark 6} - The catalogue of admissible $6$-tuples
contained in \cite{Ca} lists the complete sequence of canonical
$6$-tuples associated to prime $3$-manifolds of genus $2$, up to
complexity 21.
\end{remark}

By means of $2$-symmetric transformations, we can relate different
$H$-orbits representing the same $3$-manifold. For this purpouse,
let $\s:\FF\to\FF$ be the map $$\s\hq=\begin{cases}\hqp&\text{ if
} q_0\ne 0\\ \hq&\text{ if } q_0=0\end{cases},$$ where
$h_0',h_1',h_2',q_0',q_1',q_2'$ are the integers defined by
(\ref{4.1}), (\ref{4.2}), (\ref{4.3}) and (\ref{4.4}).

By direct computation, it is easy to check the following properties.

\begin{lemma}\label{Lemma 7} Let $\s,\z_2,\z_3:\FF\to\FF$ be the maps introduced
above. Then $\s^2=1$, $\s\z_2=\z_2\s$ and $\s\z_3=\z_3\s$.
\end{lemma}

\begin{proposition}\label{Proposition 8} Let $f$ be an admissible $6$-tuple. If $f$ is
$H$-equivalent to $f'$ then $\s(f)$ is $H'$-equivalent (and
therefore $H$-equivalent) to either $\s(f')$ or $\s(\z_1(f'))$ or
$\s(\z_1^2(f'))$.

\end{proposition}

\begin {proof} By Lemma \ref{Lemma 5} there exists $h'\in H'$ and
$e\in\{0,1,2\}$ such that $f=h'(\z_1^e(f'))$. Hence, by Lemma
\ref{Lemma 7}:
$\s(f)=\s(h'(\z_1^e(f')))=h'(\s(\z_1^e(f')))$.
\end {proof}

Let $G$ be the subgroup of $\text{Aut}(\FF)$ generated by
$\z_1,\z_2,\z_3$ and $\s$: $$G=<\z_1,\z_2,\z_3,\s>.$$ Moreover,
let $\FF_G=\FF/G$ be the set of all $G$-orbits. Each $G$-orbit is
a union of $H$-orbits and contains admissible $6$-tuples
associated with (in general) non-isomorphic graphs representing
the same $3$-manifold. Of course, if the orbits are very large, a
significant simplification in the catalogue of admissible
$6$-tuples can be achieved. The main result of the next section
supports this hope: in fact we shall prove that ``almost'' all of
$G$-orbits contain infinitely many elements.

\section {Traps and  trap-free orbits}

Let $r,s$ be positive integers such that $r\le s$; then an
admissible $6$-tuple is said to be a {\it trap of type\/} $(r,s)$
if it is $H$-equivalent to a $6$-tuple $(r,r,s;q_0,0,q_2)$ such
that:

\begin{itemize}
\item[$(*)$] $q_0+k(q_0+q_2),q_2+k(q_0+q_2)\in\{0,r+1,r+2,\ldots,s-1\}$,
for each $k\ge 0$,\end{itemize}

\noindent where both $q_0+k(q_0+q_2)$ and $q_2+k(q_0+q_2)$ are
considered mod $r+s$.\footnote{Observe that condition $(*)$ is
equivalent to the following one of ``finite type'':
\begin{itemize}\item[$(**)$] $q_0+kd,q_2+kd\in\{0,r+1,r+2,\ldots,s-1\}$,
for each $k=0,\ldots,(r+s)/d-1$, where $d=\text{GCD}(q_0+q_2,r+s)$.\end{itemize}}

\medskip

\noindent {\it EXAMPLES\/} - The $6$-tuples
$(1,1,1;0,0,0),(1,1,3;2,0,2)$ and $(1,1,2p-1;0,0,2q)$ are traps
respectively representing $\mathbf S^3,\mathbf S^1\times\mathbf S^2$ and
the lens space $L(p,q)$, for each $0<q<p$.

\medskip

\begin{lemma}\label{Lemma 9} If $f$ is a trap of type $(r,s)$, then $\s(f)$ is a
trap of the same type.
\end{lemma}

\begin {proof} $f$ is $H$-equivalent to a $6$-tuple
$f'=(r,r,s;q_0,0,q_2)$ verifying condition ($*$). By Proposition
\ref{Proposition 8}, $\s(f)$ is $H$-equivalent to either
$\s(f')=(r,r,s;-q_0,0,2q_0+q_2)$ or $\s(\z_1(f'))=\z_1(f')$ or
$\s(\z_1^2(f'))=(s,r,r;-q_2,q_0+2q_2,0)=\z_1^2(r,r,s;q_0+2q_2,0,-q_2)$.
Since $q\in\{0,r+1,r+2,\ldots,s-1\}$ if and only if
$-q\in\{0,r+1,r+2,\ldots,s-1\}$, it is easy to see that the
admissible $6$-tuples $(r,r,s;q_0+2q_2,0,-q_2)$ and
$(r,r,s;-q_0,0,2q_0+q_2)$ both verify condition $(*)$ and
therefore the statement is achieved.
\end {proof}

\begin{corollary}\label{Corollary 10} Let $f$ be an admissible $6$-tuple. Then:
\begin{itemize}
\item[a) ] If $f$ is a trap then its $G$-orbit is finite.
\item[b) ] If $f$ is not a trap then its $G$-orbit contains no trap.\end{itemize}
\end{corollary}

\begin {proof}
a) There is a finite number of traps of a fixed type.
b) Trivial.
\end {proof}

We shall call a {\it trap orbit\/} each $G$-orbit composed by
traps and a {\it trap-free orbit\/} each $G$-orbit without traps.
Here is the main result of this section.

\begin{theorem}\label{Theorem 11} Each trap-free orbit representing a $3$-manifold
of genus two contains infinitely many elements associated to
infinitely many non-isomorphic graphs.
\end{theorem}

In order to prove this theorem we define the map $\d:\FF\to\mathbf
N$, by $$\d(f)=\v(\s(f))-\v(f).$$ Note that $\d$ measures the
variation of the complexity of $f$ due to a $2$-symmetric
transformation.

From \ref{Theorem 3} we get:
\begin{equation}\label{6.1}\d\hq=\begin{cases} 0&\text{ if }
q_0=0\\ h_1-q_0&\text{ if } 0<q_0<h_0,h_2\\ q_0+h_1-h_0-h_2&\text{
if } q_0>h_0,h_2\\ h_1-h_0&\text{ if } h_0<q_0<h_2\\
h_1-h_2&\text{ if } h_2<q_0<h_0
\end{cases}.\end{equation}
Moreover, $\d$ is constant in each $H'$-orbit, since Lemma
\ref{Lemma 7} gives:
$\d(h'(f))=\v(\s(h'(f)))-\v(h'(f))=\v(h'(\s(f)))-\v(h'(f))=\v(\s(f))-\v(f)=\d(f)$,
for each $h'\in H'$.

If $f=\hq$ is an admissible $6$-tuple such that $h_0\le h_1\le
h_2$ (for example a canonical one), then we have:
\begin{equation}\label{6.2}\d(f)=\begin{cases} 0&\text{ if }
q_0=0\\ h_1-q_0&\text{ if } 0<q_0<h_0\\ h_1-h_0&\text{ if }
h_0<q_0<h_2\\ q_0+h_1-h_0-h_2&\text{ if } q_0>h_2
\end{cases};\end{equation}
\begin{equation}\label{6.3}\d(\z_1(f))=\begin{cases} 0&\text{ if } q_1=0\\
h_2-q_1&\text{ if } 0<q_1<h_0\\
h_2-h_0&\text{ if } h_0<q_1<h_1\\
q_1+h_2-h_1-h_0&\text{ if } q_1>h_1
\end{cases};\end{equation}
\begin{equation}\label{6.4}\d(\z_1^2(f))=\begin{cases} 0&\text{ if } q_2=0\\
h_0-q_2&\text{ if } 0<q_2<h_1\\
h_0-h_1&\text{ if } h_1<q_2<h_2\\
q_2+h_0-h_2-h_1&\text{ if } q_2>h_2
\end{cases}.\end{equation}

\begin{demo}
{Proof of Theorem \ref{Theorem 11}} Let $\o$ be a trap-free
orbit representing a $3$-manifold of genus two. We shall show
that, for each $\wt{f}\in\o$, there exists $\wt{f}'\in\o$ such
that $\v(\wt{f}')>\v(\wt{f})$. To achieve this fact, it suffices
to find a $6$-tuple $f'\in\o$ with the same complexity of $\wt{f}$
and such that $\d(f')>0$. In fact, in this case,
$\wt{f}'=\s(f')\in\o$ is such that
$\v(\wt{f}')=\v(f')+\d(f')>\v(f')=v(\wt{f})$.

Let $f=\hq$ be a representative of the $H$-orbit of $\wt{f}$ such
that $h_0\le h_1\le h_2$; therefore $f,\z_1(f),\z_1^2(f)$ are
$H$-equivalent to $\wt{f}$. By (\ref{6.3}) $\d(\z_1(f))>0$ whenever
$q_1\ne 0$.

Suppose now $q_1=0$, then $q_0\ne 0$ by Lemma \ref{Lemma 2}; from
(\ref{6.2}) we get $\d(f)>0$ whenever $h_0<h_1$.

It remains to examine the case $f=(h_0,h_0,h_2;q_0,0,q_2)$. Let
$T$ be the set $\{0,h_0+1,h_0+2,\ldots,h_2-1\}$. If $q_0\not\in T$
then $\d(f)>0$ and if $q_2\not\in T$ then $\d(\z_1^2(f))>0$. Let
us suppose $q_0,q_2\in T$; since $\o$ is trap-free, the set
$S=\{k>0\mid q_0+k(q_0+q_2)\not\in T \text{ or }
q_2+k(q_0+q_2)\not\in T\}$ is not empty. Let $m$ be the minimum of
$S$. Then $q_0+k(q_0+q_2)\in T$ and $q_2+k(q_0+q_2)\in T$, for
each $k=1,\ldots,m-1$, and either (a) $q_0+m(q_0+q_2)\not\in T$ or
(b) $q_2+k(q_0+q_2)\not\in T$. It is easy to check, by induction,
that
$f_m'=(\z_1\z_2\s)^m(f)=(h_0,h_0,h_2;q_0+m(q_0+q_2),0,-q_0-(m-1)(q_0+q_2))$
and
$f_m''=(\z_1\z_2\s)^m(\z_1\z_2(f))=(h_0,h_0,h_2;q_2+m(q_0+q_2),0,-q_2-(m-1)(q_0+q_2))$
(recall that $q\in T$ if and only if $-q\in T$). If (a) holds then
$\d(f_m')>0$ and if (b) holds then $\d(f_m'')>0$. This proves the
statement.\qed
\end {demo}

\begin{remark} \label{Remark 7}  - We point out that traps are really rare in the
class of admissible $6$-tuples. For example, the catalogue
enclosed in \cite{Ca} contains no traps among a list of nearly 700
canonical $6$-tuples. This shows that there are no traps of
complexity $\le 21$ representing prime $3$-manifolds of genus two.
\end{remark}

\section {Minimal $6$-tuples and roots}

The goal of producing a reduced catalogue of admissible $6$-tuple
representing all $3$-manifolds of genus two suggests looking for a
suitable representative for each $G$-orbit (a ``super-canonical''
$6$-tuple), which is possibly minimal as regards to complexity.

Let $C$ be the set of all canonical $6$-tuples. We say that $f\in
C$ is {\it minimal\/} if $\v(f')\ge\v(f)$, for each $6$-tuple $f'$
$G$-equivalent to $f$. Moreover, we say that $f\in C$ is a {\it
root\/} if $\v(f')>\v(f)$, for each $6$-tuple $f'$ $G$-equivalent
and $H$-nonequivalent to $f$.

A minimal $6$-tuple is a representative of minimal complexity of
its $G$-orbits and a root is the unique minimal $6$-tuple of the
$G$-orbit. Although not every $G$-orbit admits a root, very often
this is the case.

\begin{lemma} \label{Lemma 12} Let $f$ be a canonical $6$-tuple. Then:
\begin{itemize}
\item $\,\,f$ is minimal if and only if $\d(\z_1^i(f))\ge 0$, for $i=0,1,2$;
\item $\,\,f$ is a root if and only if $\d(\z_1^i(f))>0$ whenever
$\s(\z_1^i(f))\notin [f]_H$, for $i=0,1,2$.\end{itemize}
\end{lemma}

\begin {proof} In one direction ($\Rightarrow$) the statement is
trivial since $\s(f),\s(\z_1(f)),\s(\z_1^2(f))$ are $G$-equivalent
to $f$. To prove the converse, denote by $\S$ the graph whose
vertex-set is the set $C$ of all canonical $6$-tuples and whose
edge-set is defined by the following rule: join two different
vertices $f$ and $f'$ by an edge iff there exist two admissible
$6$-tuples $\wt{f}\in[f]_H,\wt{f}'\in[f']_H$, such that
$\wt{f}'=\s(\wt{f})$. The graph $\S$ is well-defined because
$\s^2=1$; moreover, it is an infinite graph without loops or
multiple edges. Each connected components of $\S$ corresponds to a
$G$-orbit and each vertex of $\S$ has degree $\le 3$ by
Proposition \ref{Proposition 8}: in fact, the vertices which are
adjacent to a given vertex $f$ are the canonical representatives
of the $H$-orbits $[\s(f)]_H,[\s(\z_1(f))]_H,[\s(\z_1^2(f))]_H$
distinct from $[f]_H$. We are now going to prove that if $f'$ is
adjacent to $f$ and $\v(f')<\v(f)$, then the other vertices which
are adjacent to $f$ have complexity $>\v(f)$. First of all, if
$\s(\z_1(f))$ is not $H$-equivalent to $f$ then $q_1\ne 0$ and
therefore $\d(\z_1(f))>0$ by (\ref{6.3}). Moreover, from (\ref{6.2}) we get
$\d(f)\ge 0$. Suppose now $\d(\z_1^2(f))<0$, then $h_0<h_1$ by
(\ref{6.4}) and both $\d(f),\d(\z_1(f))>0$. As a consequence, any path
in $\S$ whose sequence of vertices is $f_0=f,f_1,\ldots,f_n$ has
the following property: if $\v(f)\le\v(f_1)$ (resp.
$\v(f)<\v(f_1)$), then $\v(f_i)\le\v(f_{i+1})$ (resp.
$\v(f_i)<\v(f_{i+1})$), for each $i=0,\ldots,n-1$. Now, if
$\d(\z_1^i(f))\ge 0$ (resp. $\d(\z_1^i(f))>0$ whenever
$\s(\z_1^i(f))\notin [f]_H$), for $i=0,1,2$, then all vertices
which are adjacent to $f$ have not lower (resp. have greater)
complexity; hence, each path of positive length starting from $f$
ends in a vertex $f'$ such that $\v(f')\ge\v(f)$ (resp.
$\v(f')>\v(f)$) and therefore $f$ is minimal (resp. is a
root).
\end {proof}

As a direct consequence of Lemma \ref{Lemma 12} we can find a
complete characterization of minimal $6$-tuples and roots.

\begin{theorem}\label{Theorem 13} A canonical $6$-tuple $f=\hq$ is minimal if and
only if $$q_2<h_0\quad\text{or}\quad
q_2>h_1+h_2-h_0\quad\text{or}\quad h_0=h_1<q_2<h_2.$$ Moreover,
each minimal $6$-tuple is a root with the exception of the
following cases:

\begin{itemize}
\item[a) ] $h_0=h_1<q_2<h_2$ and $q_2\neq -q_0,(h_0+h_2)/2$ and, when $q_1=0$,
$q_2\neq (h_0+h_2)/2-q_0$;

\item[b) ] $h_0=h_1<q_0<h_2$ and $q_0\neq -q_2,(h_0+h_2)/2$ and, when $q_1=0$,
$q_0\neq (h_0+h_2)/2-q_2$.\end{itemize}

\end{theorem}

\begin {proof} From (\ref{6.2}) and (\ref{6.3}) we always get $\d(f)\ge 0$ and
$\d(\z_1(f))\ge 0$. Moreover, $\d(\z_1^2(f))\ge 0$ when either
$q_2<h_0$ or $q_2>h_2+h_1-h_0$ or $h_0=h_1<q_2<h_2$, by (\ref{6.4}).

Now, if $q_i=0$ then $\s(\z_1^i(f))\in [f]_H$, for $i=0,1,2$.
Therefore $\d(\z_1(f))>0$ whenever $\s(\z_1(f))\not\in [f]_H$. On
the other hand, it is easy to check that case a) (resp. case b))
includes all the minimal $6$-tuples $f$ such that
$\d(\z_1^2(f))=0$ and $\s(\z_1^2(f))\not\in [f]_H$ (resp. such
that $\d(f)=0$ and $\s(f)\not\in [f]_H$).
\end {proof}

%
% ---- Bibliography ----@
%

\end{article}

\end{document}